\documentstyle[12pt,leqno]{article}
\makeatletter                                            
\renewcommand{\@begintheorem}[2]{                        
\rm \trivlist \item [\hskip \labelsep {\bf #2\ \ #1.}]   
				}                        
\makeatother                                             
\makeatletter
\def\section{\@startsection {section}{1}{\z@}{-3.5ex plus -1ex minus
 -.2ex}{1.5ex plus .2ex}{\large\bf}}
\def\subsection{\@startsection{subsection}{2}{\z@}{-3.25ex plus -1ex
minus
 -.2ex}{-1em}{\normalsize\bf}}
\let\emppsubsection\subsection

\def\empsubsection[#1]#2{\emppsubsection[#1]{#2\unskip}}
\def\subsection{\secdef\empsubsection{\emppsubsection*}}

\makeatother

\newcommand{\newsubsubsection}%
{{\bf\refstepcounter{subsubsection}\thesubsubsection\ \ }}

\newcommand{\Bbb}{\bf}
\newcommand{\ts}{\vspace{\baselineskip}\noindent{\bf Proof.}$\;\;$}

\newcommand{\CC}{{\Bbb C}}

\newcommand{\PP}{{\Bbb P}}

\newcommand{\ZZ}{{\Bbb Z}}
\newcommand{\al}{\alpha}
\newcommand{\e}{\epsilon}

\newcommand{\la}{\lambda}
\newcommand{\om}{\omega}
\newcommand{\s}{\sigma}

\newcommand{\T}{\Theta}

\newcommand{\cO}{{\cal O}}

\newcommand{\qed}{{\unskip\nobreak\hfill\hbox{ $\Box$}\par}}

\title{On the Hitchin System}

\author{{\sc Bert van Geemen}\\
University of Utrecht, The Netherlands\\
{\sc Emma Previato}\\
 Boston University, USA
\thanks {Research partially supported by NSF Grant
DMS-9105221 at Boston University and DMS-9022140 at MSRI.}
}
\date{}

\begin{document}
\maketitle

\section{Introduction}

\subsection{}
What is known as the Hitchin system is a completely integrable hamiltonian
system (CIHS)
involving vector bundles over algebraic curves, identified by Hitchin in
(\cite{H1},
\cite{H2}). It was recently generalized by Faltings \cite{F}.
In this paper we only consider the case of rank-two vector bundles with
trivial determinant. In that case the Hitchin system corresponding to a curve
$C$ of genus $g$ is obtained as follows. Let
$$
{\cal M}:=\{ E\rightarrow C:\;E\;\mbox{a semi-stable rank two bundle},\;
\wedge^2E\cong \cO\;\}/\sim_{\rm S}
$$
be the moduli space of (S-equivalence classes of) semi-stable rank-two vector
bundles on $C$. Then ${\cal M}$ is a projective variety (singular if $g>2$) of
dimension $3g-3$.

The locus of stable bundles ${\cal M}^s$  is the set of  smooth
points of ${\cal M}$ for $g>2$. The cotangent space of ${\cal M}$ at a stable
bundle $E$
is :
$$
T^*_E{\cal M}=Hom_0(E,E\otimes K),\mbox{with}\quad
Hom_0(E,E\otimes K):=H^0(C,{\cal E}nd_0(E)\otimes K)
$$
where ${\cal E}nd_0(E)$ is the sheaf of endomorphisms of $E$ with trace zero
and $K$ is the canonical bundle on $C$.
A $\Phi\in Hom_0(E,E\otimes K)$ is called a Higgs field. The determinant of a
Higgs field $det(\Phi)\in Hom(\wedge^2E,\wedge^2(E\otimes K))=H^0(C,2K)$ gives
a map
$$
det:T^*_E{\cal M}=Hom_0(E,E\otimes K)\longrightarrow
H^0(C,2K),
$$
which globalizes to a map on $T^*{\cal M}^s$. Hitchin considered the map:
$$
H:T^*{\cal M}^s\longrightarrow H^0(C,2K),\quad
\Phi\mapsto det(\Phi)
$$
and  showed that it is a CIHS
in the following sense: the functions
on $T^*{\cal M}^s$
that one obtains by choosing any basis
in $H^0(C,2K)$ are a complete set of hamiltonians in involution
(with respect to the natural symplectic
structure on a cotangent bundle). Since $det$ is homogeneous of degree two in
the fibre variables ($det(t\Phi)=t^2det(\Phi))$, one can define a
(rational) projective Hitchin map
$$
\overline{H}:\PP T^*{\cal M}^s\longrightarrow\PP H^0(C,2K)=|2K|
$$
and it is in fact this map that we consider.

\subsection{}
In the first two sections of this paper we define and study a set of Higgs
fields
associated to any semi-stable bundle $E$. These results are then applied to the
$g=2$ case;
they
may also be of independent interest for studying
    moduli spaces of Higgs bundles, which are pairs $(E,\Phi)$
    as above, with certain restrictions.

To study ${\cal M}$ as a projective variety
(see \cite{NR}, \cite{B1}, \cite{B2})
one associates to any $E\in {\cal M}$
a divisor $D_E$ in the
Jacobian of $C$
using
which, questions
on
rank two bundles are rephrased in terms of line bundles (and extensions).
We exhibit a natural map
$$
\phi_E:D^{sm}_E\longrightarrow \PP Hom_0(E,E\otimes K),\qquad
\xi\mapsto \Phi_\xi
$$
(with $D_E^{sm}$ the smooth points of $D_E$)
and we are able to compute $det(\Phi_\xi)$. The result is best stated in a
diagram (the stable case of
Proposition \ref{detPhi}): For any stable bundle $E$ the following diagram
commutes:
\begin{equation}\label{diagram}
\begin{array}{rcccl}
&&\PP T_E^*{\cal M}&&\\
&\phi_E\nearrow&&\searrow\bar{H}&\\
D_E^{sm} &&&&\PP H^0(C,2K)\\
&\psi_E\searrow&&\nearrow Sq&\\
&&\PP H^0(K)&&
\end{array}
\end{equation}
here $\psi_E$ is just the Gauss map of the divisor $D_E$ inside the Jacobian
and
$Sq(\om)=\om\otimes \om$.

Thus the divisor $D_E$ (rather, its image in $\PP T_E^*{\cal M}$) plays an
important role in the study of the fibers of $\bar{H}$ over the quadratic
differentials which are squares of one forms. However, our
fiberwise approach (for the map $T^*{\cal M}\rightarrow {\cal M}$) is, in a
sense, perpendicular to Hitchin's approach which studies the fibers of
$H:T^*{\cal M}\rightarrow H^0(C,2K)$. That approach
establishes that such a fiber, over a (general) quadratic differential $\eta$
is the Prym
variety associated with a `spectral' double cover $C_\eta\rightarrow C$
defined by $\eta$ (\cite{H2},\,\cite{BNR}). It would be interesting to relate
our results to a study of the fibers of $H$.

\subsection{}
In the remaining sections we apply these results to investigate the
case $g=2$. Then the space ${\cal M}$ is isomorphic to $\PP^3$ (\cite{NR}),
so that we look for a CIHS on $T^*\PP^3$ (and also on the open subset
$T^*\CC^3$).

Using information on $\phi_E$ from $\S$\ref{phiE}
we work out the maps of the diagram (1.2.1)
in the genus two case in $\S$\ref{g=2}.
Here we encounter some classical algebraic geometry of curves of genus two
and three. It turns out that finding $H$ explicitly involves
a problem in line geometry in $\PP^3$
(a sketch of the solution in fact appears in
in J.H. Grace's article ``Line Geometry'' in
the Encyclopaedia Britannica, 1911).
We can thus make an educated guess as to what the explicit
hamiltonians
should be. A computer calculation (using the {\it Mathematica}
system) showed that our candidates actually define a CIHS. We
are not able to show that our hamiltonians define the Hitchin map, but we
can prove that Hitchin's hamiltonians and ours differ by multiplication by
functions from the base (an open set in $\PP^3$). For a more precise result we
would have to extend the results of section \ref{g=2} to enlarge the open set
in the base were those results hold, or we would need further information on
Hitchin's
system.

\subsection{Acknowledgements.} The first named author
 wishes to thank: the University of Pavia for a six month stay where much of
the work on this paper was done, the NSF for supporting visits to Boston
University under Grant DMS-9105221, and J. de Jong for helpful discussions.

 The second named author wishes to acknowledge:
Carolyn Gordon's invitation to MSRI for two
weeks of the special year in Differential Geometry 1993/94
(research at MSRI supported in part by NSF grant \# DMS 9022140), and
participation in the LMS/Europroj Workshop ``Vector bundles in
algebraic geometry'' (Durham, 1993; organizers N. Hitchin, P. Newstead
and W.M. Oxbury), on which occasion N. Hitchin provided generous
insight.

\section{Higgs fields}

\subsection{}
We fix some notations and recall some basic facts.

In this paper $C$ will be a smooth, irreducible projective curve
of genus $g>1$ over $\CC$
and $E$ will be a rank two semi-stable bundle on $C$ with trivial
determinant.

Since $\wedge^2E\cong {\cal O}$, we have:
$$
E\wedge E={\cal O}_C,\qquad {\rm so}\quad
E\cong E^*:={\cal H}om(E,{\cal O}_C),\quad e\mapsto [f\mapsto e\wedge f]
$$
thus $E$ is self-dual.
This gives isomorphisms:
$$
E\otimes E\cong E^*\otimes E={\cal E}nd(E),\qquad
S^2E\cong {\cal E}nd_0(E)
$$
with ${\cal E}nd_0(E)\subset {\cal E}nd(E)$  the sheaf of endomorphisms of
trace zero.
We recall that
$End_0(E):=H^0(C,{\cal E}nd_0(E))=0$ for a stable bundle $E$, the only
endomorphisms of $E$ being
scalar multiples of the identity. Thus:
$ H^0(C,S^2E)=0$.

For a vector space $V$, we let $\PP V$ be the space of one dimensional
linear subspaces of $V$.

\subsection{}\label{inv}
We will construct Higgs fields by relating $E$ to line bundles. Such a
connection is provided by the following results.
Let $E$ be a semi-stable rank two bundle on $C$ with $det(E)=\cO$.
Associated to $E$ is a divisor (\cite{B1}, 2.2):
$$
D_E:=\{\xi\in Pic^{g-1}(C):\;\dim H^0(\xi\otimes E)>0\;\}.
$$
With its natural scheme structure, $D_E$ is linearly equivalent to $2\Theta$.
Here $\Theta$ is the natural theta divisor:
$$
\Theta:=\{\xi\in Pic^{g-1}(C):\;\dim H^0(C,\xi)>0\;\},\qquad
D_E\in |2\T|.
$$

On $Pic^{g-1}(C)$ there is a natural involution:
$$
\iota:Pic^{g-1}(C)\longrightarrow Pic^{g-1}(C),\qquad
\xi\mapsto  K\otimes
\xi^{-1}.
$$

All divisors in $|2\T|$ are invariant under the involution $\iota$. For a $D_E$
that is easy to check since by Riemann-Roch and Serre duality:
$$
\dim H^0(\xi\otimes E)=\dim H^1(\xi\otimes E)=
\dim H^0(\xi^{-1}\otimes K\otimes E).
$$
Note that $H^0(\xi^{-1}\otimes K\otimes E)=Hom(\xi,E\otimes K)$ and that, since
$E$ is self-dual,
$H^0(\xi\otimes E)=Hom(E,\xi)$.
Thus we have:
$$
\xi\in D_E\Longleftrightarrow Hom(E,\xi)\neq 0 \Longleftrightarrow
Hom(\xi,E\otimes K)\neq 0.
$$

\subsection{}
For any semi-stable $E$ (cf.\ \cite{L}, Cor. V.6):
$$
\xi\in D^{sm}_E\;\Longrightarrow\; \dim H^0(\xi\otimes E)=1.
$$
Thus for $\xi\in D^{sm}_E$ there
are unique (up to scalar multiple) maps:
$$
\pi:E\longrightarrow \xi,\qquad
\tau: \xi\longrightarrow E\otimes K.
$$
The composition
$$
\tau\circ \pi:\;E\longrightarrow E\otimes K
$$
is an element of $Hom(E,E\otimes K)$, defined (up to scalar multiple) by
$\xi$.

\subsection{Definitions.} \label{defpp}
Let $E$ be a semi-stable rank-two bundle on $C$ with $det(E)=\cO$.
We define rational maps:
$$
\phi_E:\; D_E^{sm}\longrightarrow \PP Hom_0(E,E\otimes K),
$$
$$
\xi\mapsto \Phi_\xi:=
\tau\circ \pi-(1/2)(id_E\otimes tr(\tau\circ \pi)):
E\otimes\cO\longrightarrow E\otimes K
$$
and
$$
\psi_E:\; D_E^{sm}\longrightarrow \PP H^0(C,K)\;=\PP Hom(\xi,\xi\otimes K),
$$
$$
\psi_E(\xi)=(\pi\otimes id_K)\circ\tau:\;
\xi\stackrel{\tau}{\longrightarrow} E\otimes K
\stackrel{\pi\otimes 1}{\longrightarrow}\xi\otimes K .
$$

\subsection{} Now we have a large supply of Higgs fields, the $\Phi_\xi$'s.
It is surprisingly easy to determine $\psi_E$. In Proposition \ref{detPhi}
we will see how that already determines the Hitchin map to a large extent.

Recall that the cotangent bundle to $Pic^{g-1}(C)$ is trivial:
$$
T^*Pic^{g-1}(C)\cong Pic^{g-1}(C)\times H^1(Pic^0(C),\cO)^*.
$$
 For a smooth point $\xi$ in a divisor $D\subset Pic^{g-1}(C)$ the tangent
space
to $D$ at $\xi$ is then defined by an element of $H^1(Pic^0(C),\cO)^*$, unique
up to scalar multiple. The corresponding morphism $D^{sm}\rightarrow \PP
H^1(Pic^0(C),\cO)^*$ is called the Gauss map.

\subsection{Proposition.}\label{psig}
 The map $\psi_E$ is  the Gauss map on $D_E\subset Pic^{g-1}(C)$.
$$
\psi_E:D_E^{sm}\longrightarrow \PP H^1(Pic^0(C),\cO)^*=\PP H^0(C,K).
 $$
In particular, $\psi_E$ is a morphism.

\ts
By \cite{L}, Prop.\ V.2, we know that for $\xi\in D_E^{sm}$ the space
$$
T_\xi D_E\subset T_\xi Pic^{g-1}(C)= H^1(Pic^0(C),\cO)=H^1(C,\cO)
$$
is defined by the image of the cup-product map
$$
H^0(C,E\otimes\xi)\otimes H^0(C,E\otimes \iota(\xi))\longrightarrow
H^0(C,K)\cong H^1(C,\cO)^*.
$$
This map coincides with the composition:
$$
Hom(E\otimes  K,\xi\otimes K)\otimes Hom(\xi,E\otimes K)
\longrightarrow Hom(\xi,\xi\otimes K)\cong H^0(K),
$$
and in our case we recover the definition of $\psi_E$:
$$
(\pi\otimes 1)\otimes \tau\mapsto (\pi\otimes 1)\circ\tau =\psi_E(\xi).
$$
(One may in fact also consider $Hom(E,\xi_\eta)$ where $\xi_\eta$ is a
deformation of $\xi$ given by $\eta\in H^1(C,\cO)$. Then $\eta\in T_\xi D_E$
iff
$\pi\in Hom(E,\xi)$ lifts to $Hom(E,\xi_\eta)$ iff $\pi\cup \eta=0\in
Ext^1(E,\xi)=H^1(E\otimes\xi)$, which gives the statement above. The
justification for this argument is given in \cite{L}, II.)
\qed

\subsection{} We are interested in computing the determinant of the Higgs
field $\Phi_\xi$. Since the maps in \ref{defpp} are only defined up to scalar
multiple,
we consider
$$
det: Hom_0(E,E\otimes K) \longrightarrow \PP H^0(C, 2K),
\qquad \Phi\mapsto \langle det(\Phi)\rangle.
$$
Let
$$
Sq: H^0(C,K)\longrightarrow H^0(C, 2K),\qquad
\om\mapsto\om^{\otimes 2}.
$$

\subsection{Proposition.}\label{detPhi}
 For a semi-stable $E$ and $\xi\in D_E^{sm}$ we have:
$$
det(\Phi_\xi)=\psi_E(\xi)^{\otimes 2}\qquad
(\in \PP H^0(C, 2K)).
$$
Thus, the compositions $det\circ\phi_E$ and $Sq\circ\psi_E$ coincide:
$$
\begin{array}{cccccc}
det\circ \phi_E: & D_E^{sm}&\stackrel{\phi_E}{\longrightarrow}&
\PP Hom_0(E,E\otimes K)&
\stackrel{det}{\longrightarrow}& \PP H^0(C, 2K),\\
Sq\circ\psi_E:& D_E^{sm}&\stackrel{\psi_E}{\longrightarrow} &\PP H^0(C,K)&
\stackrel{Sq}{\longrightarrow} &\PP H^0(C, 2K).
\end{array}
$$

\ts
Since $\psi_E:D_E^{sm}\rightarrow \PP H^0(C,K)$ is a morphism,
 $\psi_E(\xi)$ is (represented by) a non-zero differential form
for each $\xi\in D_E^{sm}$.
Define a canonical divisor on $C$ by:
$$
K_\xi:=(\psi_E(\xi)),\qquad {\rm let}\quad U:=C-Support(K_\xi),
$$
On the open set $U$, the map (defined by)
$\psi_E(\xi):\xi\rightarrow \xi\otimes K$ is an isomorphism. Its
inverse, composed with $\tau:\xi\rightarrow E\otimes K$, gives a
map $\xi\otimes K\rightarrow E\otimes K$ which splits the map
$\pi\otimes1:E\otimes K\rightarrow \xi\otimes K$.
$$
\begin{array}{ccc}
\xi                             & {}_{\psi_E(\xi)}& \\
\tau\Big\downarrow\phantom{\tau}&\searrow
&\\
E\otimes K          &\stackrel{\pi\otimes 1}{\longrightarrow}& \xi\otimes K
\end{array}
$$

Thus over $U$, the bundle $E$ splits:
$$
{E}_{|U}\cong L\oplus\xi_{|U},\qquad{\rm with}\quad
L:=\ker(\pi_U:E_{|U}\longrightarrow \xi_{|U})
$$
and $L$ is a line bundle on $U$.
Since $\Phi_\xi:=\tau\pi-(1/2)tr(\tau\pi)$, we get
$\Phi_{\xi|U}(L)=0$ so that:
$$
\Phi_{\xi|U}=\left(\begin{array}{cc}
-(1/2)\psi_E(\xi)&\ast\\0&(1/2)\psi_E(\xi)
\end{array}\right).
$$
Then $det(\Phi_{\xi{|U}})=-(1/4)\psi_E(\xi)^{\otimes 2}$, which
does not
vanish at any point of $U$.

If $D_E$ is irreducible,
the image of the Gauss map on $D^{sm}_E$ contains an open subset of
$\PP H^0(C,K)$.
Thus, for general $\xi$ on such a $D_E$,
$2K_\xi$ is the only divisor in $\PP H^0(C,2K)$ with support in $C-U$.
Since $det(\Phi_\xi)\in H^0(C,2K)$ and since its divisor must have support in
$C-U$, we conclude that
$$
(det(\Phi_\xi))=2K_\xi.
$$

 For general $C$, the rank of the N\' eron-Severi
group of $Pic^{g-1}(C)$ is one,
and then the reducible divisors in $|2\T|$ are a subvariety of dimension $g$
(they are unions of two translates of $\T$).
Since $\Delta({\cal M})$ (see \ref{delta}) has dimension $3g-3$ (\cite{B1}),
 the divisor $D_E$ is irreducible for general $C$.
Thus if we work in a family
of $D_E$'s over a general family of curves containing the given $D_E$,
the maps $\xi\mapsto det(\Phi_\xi)$ and $\xi\mapsto \psi_E(\xi)^{\otimes
2}$ agree on a non-empty open subset, hence must agree everywhere.
\qed

\section{The map $\phi_E$} \label{phiE}

\subsection{}\label{delta}
In this section we consider only stable bundles $E$ and we
 study the map $\phi_E:D_E^{sm}\rightarrow \PP T^*_E{\cal M}$.
We use the codifferential of the map:
$$
\Delta: {\cal M}\longrightarrow |2\T|,\qquad E\mapsto D_E
$$
to relate the cotangent bundles of ${\cal M}$ and of the projective space
$|2\T|$.
 First of all we recall some facts on the cotangent bundle to a projective
 space
and on the dual of $|2\T|$ (following \cite{NR2}, $\S$ 3).

\subsection{} \label{dual}
Let $V$ be a vector space.
The dual of the Euler sequence on $\PP V$ gives:
$$
0\longrightarrow T^*\PP V\longrightarrow V^*\otimes{\cal O}(-1)\longrightarrow
{\cal O}
\longrightarrow 0,
$$
the last non-trivial map is given by
$(\ldots,s_i,\ldots)\mapsto\ldots+x_is_i+\ldots$ over
$(\ldots :x_i:\ldots)\in\PP V$.

Taking the associated projective bundles
we have an isomorphism:
$$
\PP T^*\PP V\cong I:=\{\;(x,h)\in \PP V\times \PP V^*:\;
x\in h\;\},
$$
the variety $I$ is called the incidence bundle.

\subsection{}\label{delta2}
To identify the dual of $|2\T|$ we use the map:
$$
\delta:Pic^{g-1}(C)\longrightarrow |2\T_0|,\qquad
\xi\mapsto D_\xi:=L^*_\xi\Theta +L^*_{\iota(\xi)}\T,
$$
here $\Theta_0$ is (any) symmetric theta divisor in $Pic^0(C)$ (the linear
equivalence class of $2\Theta_0$ is independent of the choice) and
$$
L_\al:Pic(C)\longrightarrow Pic(C),\qquad
\beta\mapsto \al\otimes\beta
$$
is translation by $\al$ in $Pic(C)$.

Pulling back the linear forms on $|2\Theta_0|$ gives an isomorphism
$$
\delta^*:H^0(|2\Theta_0|,{\cal O}(1))=H^0(Pic^0(C),2\Theta_0)^*
\stackrel{\cong}{\longrightarrow} H^0(Pic^{g-1}(C),2\Theta).
$$
Projectivizing gives
$\delta^*:|2\Theta_0|^*\stackrel{\cong}{\longrightarrow} |2\Theta|$.
In fact, there is a commutative diagram:
$$
\begin{array}{ccc}
& & |2\T|^*\\
&{}^\nu\nearrow\phantom{{}^\nu}& \\
Pic^{g-1}(C)& &\;\;\downarrow (\delta^*)^*\\
& {}_\delta\searrow\phantom{{}_\delta}&\\
& &|2\T_0|
\end{array}
$$
where $\nu$ is the natural map (see \cite{NR2} and \cite{B2} \S 2 for a
variant). From now on, $I$ will be the incidence bundle:
$$
I:=\PP T^*|2\T|\subset |2\T|\times |2\T_0|.
$$
We will denote by
 $(d\Delta)^*$ the projectivized codifferential of $\Delta$
(see \ref{delta}):
$$
(d\Delta)^*:\PP T^*|2\T|= I\longrightarrow \PP T^*{\cal M}.
$$

\subsection{Lemma.}\label{dd}
 Let $D\in |2\T|$ and let $\xi\in Pic^{g-1}(C)$. Then
$$
(D,D_\xi)\in I\;\Longleftrightarrow\; \xi\in D.
$$

\vspace{\baselineskip}
\noindent
{\bf Proof.}$\;\;$
 For $\xi\in Pic^{g-1}(C)$, $\nu(\xi)\in |2\T|^*$
is the hyperplane in
$|2\T|$ consisting of the
divisors passing through $\xi$.
Thus $D\in |2\T|$ and $\xi\in |2\T|^*$ are incident iff $\xi\in D$.
The dual of the isomorphism $\delta^*$ maps $\nu(\xi)$ to $D_\xi$ so the result
follows.
\qed

\subsection{} We recall that for non-hyperelliptic curves,
the map $\Delta$ has degree one over its image \cite{B1}
(so it is locally an isomorphism with its image for generic $E$) and is an
embedding for the general curve as recently announced by Y. Laszlo
and also by
S. Brivio and A. Verra jointly.
In case $g=2$ the map is an isomorphism \cite{NR} but for hyperelliptic curves
of genus greater than two
the map is 2:1 and `ramifies' along a subvariety of dimension $2g-1$
\cite{B1}.

\subsection{Proposition.}\label{glue}
Let $E$ be a stable bundle such that the map $\Delta$ is locally at $E$ an
isomorphism with its image. Then
the rational map
$$
\phi_E:D_E\longrightarrow \PP T^*_E{\cal M},\qquad \xi\mapsto \Phi_{\xi},
$$
is the left-hand column in the diagram:
$$
\begin{array}{ccc}
D_E&\hookrightarrow&Pic^{g-1}(C)\\
\Bigg\downarrow&&\Bigg\downarrow\delta\\
\PP T^*_{D_E}|2\T| &\hookrightarrow& |2\T_0|\\
(d\Delta_E)^*\Bigg\downarrow\phantom{(d\Delta_E)^*} & &\\
\PP T^*_E {\cal M}&&\\
\end{array}
$$
where the last vertical arrow is a linear projection given by the dual of the
differential of $\Delta$ at $E\in{\cal M}$.

\vspace{\baselineskip}
\noindent
{\bf Proof.}$\;\;$
Let $\xi\in D_E$. Then $\delta(\xi)=D_\xi\in|2\T_0|=|2\T|^*$ corresponds to a
hyperplane $H_\xi\subset |2\T|$. By Lemma \ref{dd}, $\xi\in D_E$ implies
$D_E\in H_\xi\subset |2\T|$.
This says that $H_\xi$ passes through
$D_E=\Delta(E)\in\Delta({\cal M})\subset |2\T|$ and
(by the assumption on  local isomorphism)
 defines a codimension
$\leq 1$ subspace in $T_E{\cal M}$. We must show that
$\Phi_\xi\;(\in T^*_E{\cal M})$ is the defining equation for this subspace.

We first determine $H_\xi\cap \Delta({\cal M})$;
the pull-back $\Delta^*H_\xi$ will be the divisor
$\tilde{D}_\xi$ defined below; in particular, the subspace of $T_E{\cal M}$
defined by
$H_\xi$ is $T_E\tilde{D}_\xi$.

We recall from \cite{B1} that
$$
Pic({\cal M})\cong\ZZ,\qquad{\rm and}\quad {\cal L}:=\Delta^*({\cal O}(1))
$$
is the ample generator of this group. Moreover, the natural map
$$
{\cal M}\longrightarrow \PP H^0({\cal M},{\cal L})^*
$$
actually coincides with $\Delta$.

Define for $\xi\in Pic^{g-1}(C)$:
$$
\tilde{D}_\xi:=\{E\in{\cal M}:\; H^0(C,E\otimes\xi)\neq 0\;\}.
$$
This divisor, with its natural scheme structure, is defined by a section of
$H^0({\cal M},{\cal L})$.

Restriction to the Kummer variety of $Pic^0(C)$ (= locus of non-stable bundles)
in ${\cal M}$ induces the isomorphism (\cite{B1}):
$$
 \PP H^0({\cal M},{\cal L})\stackrel{\cong}{\longrightarrow} |2\T_0|,
\qquad {\rm and}\quad \tilde{D}_\xi\mapsto D_\xi
$$
(indeed, for $L\in Pic^0(C)$ one has $H^0((L\oplus L^{-1})\otimes\xi)>0$ iff
$L\in L_\xi^*\T$ or $L\in L_{\iota(\xi)}^*\T$ iff $L\in D_\xi$).
Now, by definition, the hyperplane $H_\xi$ intersects $Pic^0(C)$ in $D_\xi$ so
that
 $H_\xi$ intersects ${\cal M}$ in $\tilde{D}_\xi$, as desired.

We must now show that $\Phi_\xi$ defines the subspace $T_E \tilde{D}_\xi$ of
$T_E{\cal M}$. The divisor $\tilde{D}_\xi$ is (the closure of) the image of the
(rational) map:
$$
\rho:\PP H^1(\xi^{-2})\longrightarrow {\cal M},\qquad
\e\mapsto [E_\e]
$$
where $E_\e$ is the extension defined by $\e\in
Ext^1(\xi,\xi^{-1})=H^1(\xi^{-2})$,
\begin{equation}
\label{ext}
0\longrightarrow \xi^{-1}\longrightarrow E_\e\stackrel{\pi}{\longrightarrow}
\xi
\longrightarrow 0.
\end{equation}
These maps were studied in detail by Bertram in \cite{B}.
We will now assume $\pi:E\rightarrow \xi$ to be surjective, so $\xi^{-1}$ is a
subbundle  of $E$. By specialization the result
follows for all $\xi\in D_E$, all $E$.

We tensor the sequence \ref{ext} by $E$,
obtaining the following sequence for $S^2E$:
$$
0\longrightarrow \xi^{-2}\longrightarrow
S^2E \longrightarrow \xi\otimes E
\longrightarrow 0.
$$
Since ${\cal E}nd_0(E)\cong S^2E$,
the differential of $\rho$ is the natural map:
$$
d\rho:H^1(\xi^{-2})/H^0(\xi\otimes E)\longrightarrow H^1(S^2E)=T_E{\cal M},
\qquad
{\rm so}\quad
T_E \tilde{D}_\xi=Im(d\rho).
$$

Dualizing the sequence above and tensoring it by $K$ we get:
$$
0\longrightarrow \xi^{-1}\otimes E\otimes K\longrightarrow
S^2E\otimes K \longrightarrow \xi^{ 2}\otimes K
\longrightarrow 0.
$$
Now we have:
$$
\begin{array}{rcl}
\langle\Phi_\xi\rangle
&=& Im(H^0(\xi^{-1}\otimes E\otimes K)=\langle\tau\rangle
	  \longrightarrow H^0(S^2E\otimes K))\\
&=& \ker(H^1(E\otimes\xi)\longleftarrow H^1(S^2E))^*\\
&=& Im(H^1(S^2E)\longleftarrow H^1(\xi^{-2}))^*.
\end{array}
$$
Therefore we have indeed, as claimed:
$$
T_E \tilde{D}_\xi=\ker(\Phi_\xi).
$$
\qed


\subsection{}\label{covD}
The previous proposition shows that the line bundle
$\cO_{Pic^{g-1}(C)}(2\T)_{|D_E}\cong \cO_{D_E}(D_E)$ plays an essential role as
regards
the map $\phi_E$. Recall that the divisor $D_E$ is invariant under
the involution $\iota$ (see \ref{inv}).
The following lemma shows how $\iota$ acts on the global sections of this line
bundle.

\subsection{Lemma.}\label{split}
The  involution $\iota$ gives a splitting in an invariant and an
anti-invariant part:
$$
\begin{array}{rcccc}
H^0(D_E,\cO_{D_E}(D_E))&=& H^0(D_E,\cO_{D_E}(D_E))_+&\oplus&
H^0(D_E,\cO_{D_E}(D_E))_-\\
&\cong&H^0(Pic^{g-1}(C),\cO(D_E))/\langle s_E\rangle & \oplus&
H^1(Pic^{g-1}(C),\cO).
\end{array}
$$
Here $s_E\in H^0(Pic^{g-1}(C),\cO(D_E))$ is a section with divisor
$(s_E)=D_E$.

The map $\phi_E$ factors over the natural map
$$
D_E\longrightarrow \PP H^0(D_E,\cO_{D_E}(D_E))_+^*\;(\cong \PP
T^*_{D_E}|2\T|).
$$
The map $\psi_E$ (the Gauss map)
is the natural rational map:
$$
D_E\longrightarrow \PP H^0(D_E,\cO_{D_E}(D_E))_-^*\cong
 \PP H^1(Pic^{g-1}(C),\cO)^*.
$$
Therefore both $\phi_E$ and $\psi_E$ factor over $\bar{D}_E$.

\ts
The exact sequence of sheaves on $Pic^{g-1}(C)$:
$$
0\longrightarrow \cO\stackrel{\cdot s_E}{\longrightarrow}
\cO(D_E)\longrightarrow \cO_{D_E}(D_E)\longrightarrow 0
$$
gives the cohomology sequence:
$$
0\longrightarrow H^0(Pic^{g-1}(C),\cO(D_E))/\langle s_E \rangle
\longrightarrow H^0(D_E,\cO_{D_E}(D_E))
\longrightarrow H^1(Pic^{g-1}(C),\cO)\longrightarrow 0.
$$
It is well known that $\iota^*$ acts as the identity on
$H^0(Pic^{g-1}(C),\cO(D_E))$ and as minus the identity on
$H^1(Pic^{g-1}(C),\cO)$. The remaining assertions are standard.
\qed

\section{The genus two case}\label{g=2}

\subsection{} To determine the Hitchin map in the genus two case, we study
first
the divisors $D_E$ for general stable $E$ and we
study three of the four maps from the diagram \ref{diagram} (see
\ref{lower2}, \ref{upper}).
This leads to quite classical geometry involving for instance \' etale double
covers and tangent
conics. Then we can easily determine the fourth map, which is Hitchin's map
(projectivized and restricted to $\PP T^*_E{\cal M}$).
We then `rigidify' our construction using the classical Proposition
\ref{prym}.

  From now on, $C$ will be a genus two curve.

\subsection{}
The main result of \cite{NR} is that the map
$$
\Delta:{\cal M}\stackrel{\cong}{\longrightarrow} |2\Theta|\cong\PP^3,\qquad
E\mapsto D_E
$$
is an isomorphism, so ${\cal M}\cong \PP^3$.

In particular, any element in $|2\T|$ is a $D_E$ for some $E$.
Since this linear system is base-point free, the divisors,
now in fact curves,
$D_E$ are smooth and have genus 5
for general stable $E$
(a description of the singular curves in $|2\Theta|$ can be found in
\cite{Verra}).

\subsection{}\label{cov}
Each divisor in $|2\T|$ is fixed by $\iota$ and
for general $E$,
the involution $\iota$ restricted to $D_E$ is a fixed-point free involution on
a smooth curve. The induced covering
$$
\pi_E:D_E\longrightarrow \bar{D}_E:=D_E/\iota
$$
is an \' etale 2:1 covering (the associated Prym variety is $Jac(C)$,
see for example \cite{Verra}, p.\ 438).
In particular, for general $E$, $\bar{D}_E$ is a smooth genus three curve.
   From now on we will consider only such $E$.

\subsection{}\label{bunde}
The kernel of the map $\pi_E^*:Pic(\bar{D}_E)\rightarrow Pic(D_E)$ is generated
by a point $\alpha$ of order two.
One has:
$$
\pi_E^* K_{\bar{D}_E}\cong K_{D_E},\quad
\pi_{E*}\cO_{D_E}\cong\cO_{\bar{D}_E}\oplus\alpha.
$$
The adjunction formula on $Pic^{g-1}(C)$ shows $\cO_{D_E}(D_E)\cong
K_{D_E}$.
The involution $\iota$ gives a splitting in an invariant and an
anti-invariant part:
$$
H^0(D_E,K_{D_E})\cong H^0(\bar{D}_E,\pi_* K_{D_E})\cong\;
H^0(\bar{D}_E,K_{\bar{D}_E})\oplus
H^0(\bar{D}_E,K_{\bar{D}_E}\otimes\alpha),
$$
(projection formula)
which coincides with the splitting given in Lemma \ref{split}. In particular,
since $H^1(Pic^{g-1}(C),\cO)\cong H^1(C,\cO)\cong H^0(C,K)^*$, we have a
natural identification
$$
H^0(C,K)\cong H^0(\bar{D}_E,K_{\bar{D}_E}\otimes\alpha)^*.
$$

\subsection{} \label{lower2}
We will now write $C_3$ for $\bar{D}_E$, $K_3$ for $K_{\bar{D}_E}$.
The Gauss map on $D_E$ then factors over $C_3$, and on $C_3$ coincides with the
natural map
$$
C_3\longrightarrow \PP H^0(C_3,K_3\otimes\al)^*\cong
\PP H^1(Pic^{g-1}(C),\cO)^* \cong \PP^1
$$
which is therefore also essentially $\psi_E$.

The map $Sq:|K|\rightarrow |2K|$ from diagram \ref{diagram}
corresponds to the second Veronese map which embeds $\PP^1$ as a conic in
$\PP^2$. The (three dimensional) space $S^2H^0(C_3,K_3\otimes\al)^*$ may be
identified with a quotient of the (six dimensional) $H^0(C_3,2K_3)^*$
(note $2(K_3\otimes\alpha)\equiv 2K_3$), thus we have a diagram (where the last
map is a linear projection):
$$
\begin{array}{ccccccccc}
D_E& &\stackrel{\psi_E}{\longrightarrow}&
&|K|&\stackrel{Sq}{\longrightarrow}&|2K|&   &     \\
\Big|\,\!\Big|&& && \;\;\Big\downarrow \cong && \;\Big\downarrow\cong& &  \\
D_E&\stackrel{\pi_E}{\longrightarrow}\!&C_3&\!\rightarrow&\!\PP
H^0(C_3,K_3\otimes\al)^*&\longrightarrow&\!\PP S^2H^0(C_3,K_3\otimes\al)^*&
\leftarrow &\!\PP H^0(C_3,2K_3)^*.
\end{array}
$$

\subsection{}\label{upper}
As $\Delta:
{\cal M}{\rightarrow} |2\T|=\PP^3$,
is an isomorphism, the cotangent bundle to ${\cal M}$ is the incidence
bundle. The map $(d\Delta)^*$ induces an isomorphism.
$$
(d\Delta)^*:I=\{(x,h)\in \PP^3\times(\PP^3)^*:\;x\in h\;\}
\stackrel{\cong}{\longrightarrow} \PP T^*{\cal M} \quad
{\rm and}\quad \PP T_{D_E}^*|2\T|\cong \PP T^*_E{\cal M}.
$$
 From Lemma \ref{split} and \ref{bunde} we get:
$$
\phi_E:D_E\stackrel{\pi_E}{\longrightarrow} C_3
\stackrel{\kappa}{\longrightarrow} \PP H^0(C_3, K_3)^*\cong \PP T_E^*{\cal M}
$$
where $\kappa$ is just the canonical map.

\subsection{}\label{kum}
We show how the various $C_3$'s ($=\bar{D}_E$'s) fit together as $E$ moves
over $\PP^3\;(={\cal M})$.
The image $S$ of the map
$$
\delta:Pic^{g-1}(C)\longrightarrow S\subset |2\T_0|
= \PP H^0(Pic^{g-1}(C),2\T)^*=(\PP^3)^*.
$$
(see  \ref{delta2})
is the Kummer surface of $Pic^{g-1}(C)$:
$$
S\cong Pic^{g-1}(C)/\iota,\qquad \iota:L\mapsto K\otimes L^{-1}.
$$
The surface $S$ is a quartic surface and its singular locus consists of the 16
fixed points of $\iota$ (which are the theta characteristics on $C$).

Moreover, we can view $\PP T_E^*{\cal M}$ as a plane in $(\PP^3)^*$:
$$
\PP T_E^*{\cal M}=\{h\in(\PP^3)^*:\; E\in h\;\}.
$$

Proposition \ref{glue} shows:
$$
\PP T_E^*{\cal M}\cap S\,=\; \phi_E(D_E),\qquad
({\rm with}\quad \phi_E:D_E\longrightarrow \PP T_E^*{\cal M}).
$$

Thus $\phi_E(D_E)$ is a hyperplane section
of the Kummer surface $S$, hence a quartic plane curve.
For general $E$, this curve will be smooth (i.e.\ the curve $\bar{D}_E$ is
non-hyperelliptic).
 We consider only these $E$.

\subsection{}\label{hf3}
The curve $C_3$ is now non-hyperelliptic by assumption, so
 the canonical map $\kappa$
is an embedding, and the image of $C_3$ is a smooth quartic in $\PP^2$.
Pull-back along $\kappa$ gives an isomorphism
$H^0(C_3,2K_3)\cong H^0(\PP^2,\cO(2))$.

Let $s,\,t$ be a basis of $H^0(C_3,K_3\otimes\al)$. We define
conics $Q_i$ in $\PP^2$
by:
$$
s\otimes s=Q_1,\quad s\otimes t=Q_2,\quad t\otimes t=Q_3
$$
and the essential part (that is, on $D_E/\iota=C_3$) of $Sq\circ\psi_E$ is
now:
$$
C_3\longrightarrow \PP S^2H^0(C_3,K_3\otimes\al)^*\cong \PP^2,
$$
$$
x\mapsto (s^2(x):st(x):t^2(x))=(Q_1(\kappa(x)):
Q_2(\kappa(x)):Q_3(\kappa(x))).
$$
Since $Sq\circ \psi_E=\bar{H}\circ\phi_E$,
we conclude that the Hitchin map is given by:
$$
\bar{H}:\PP T_E^*{\cal M}
\longrightarrow \PP H^0(2K)=|2K|,\qquad
p\mapsto (Q_1(p):Q_2(p):Q_3(p))
$$
(since  $\kappa(C_3)=\phi_E(D_E)$ has degree 4, spans $\PP T_E^*{\cal M}$
and $\bar{H}$ has quadratic coordinate functions, $\bar{H}$
is determined by its restriction to $\phi_E(D_E)$).
Note that the inverse image in $\PP T_E^*{\cal M}$
of the conic $Sq|K|\subset |2K|$ under $\bar{H}$ is a quartic curve containing
$\phi_E(D_E)$ and thus is equal to $\phi_E(D_E)$.

\subsection{}\label{bit}
We study the construction above a little more closely and exhibit natural (up
to
scalar multiple) subsets of $H^0(C_3,K_3\otimes\alpha)^*$ and $H^0(K)$ which
correspond under the isomorphism $H^0(C_3,K_3\otimes\alpha)^*\cong H^ 0(K)$
that we found in \ref{bunde}.

For any $a,\,b\in\CC$ we have a section $as+bt\in H^0(C_3,K_3\otimes\alpha)$,
let:
$$
Q_{(a:b)}:=S^2(as+bt)=a^2Q_1+2abQ_2+b^2Q_3.
$$
The $Q_{(a:b)}$ are a quadratic system of conics, each of which is
tangent to $C_3$ (that is, has even intersection multiplicity at each
intersection point) because $Q_{(a:b)}$ cuts out twice the divisor of the
section $as+bt\in H^0(K_3\otimes \alpha)$.

There are 6 conics $H_i$, $i\in\{1,\ldots ,6\}$, in the quadratic system of
tangent conics which
split as pairs of bitangents. They correspond to the six points
$$
\langle s_i\rangle=\langle a_is+b_it\rangle\in\PP^1=
\PP H^0(K_3\otimes\al),\qquad{\rm with}\quad
det(a_i^2Q_1+2a_ib_iQ_2+b_i^2Q_3)=0
$$
(where we now view the $Q_i$ as $3\times 3$ matrices).
In this way obtain 12 bitangents of $C_3$. The other 16 bitangents are best
seen
by identifying $C_3$ with a hyperplane section of the Kummer surface $S$
(\ref{kum}).
In fact the divisor $\T$ and its translates by points of order two map
to conics in $S$, the plane through such a conic intersects $S$ in a double
conic, and thus intersects the plane in
which $\bar{D}_E$ lies in a bitangent of $\bar{D}_E$.

The following classical result relates these six points
$\langle s_i\rangle$ to the Weierstrass points of the curve $C$.

\subsection{Proposition.}\label{prym}
Under the natural isomorphism from \ref{bunde}:
$$\
\PP H^0(C_3,K_3\otimes\alpha)\stackrel{\cong}{\longrightarrow} \PP H^0(K)^*,
$$
the six points which correspond to pairs of bitangents are mapped to the six
linear maps corresponding to the Weierstrass points $p_i$
of $C$:
$$
\langle s_i\rangle\mapsto \langle [\om\mapsto \om(p_i)]\rangle
\qquad\qquad (\om\in H^0(K)).
$$

\ts
We recall the way $D_E$ and $C$ can be recovered from
$H^0(C_3,K_3\otimes\alpha)$.
It will be shown that the double cover of $\PP H^0(C_3,K_3\otimes\alpha)$
branched over the 6 points corresponding to pairs of bitangents is
isomorphic to $C$, which proves the Proposition.

Since we have $(st)^{ 2}=(s^ {2})( t^{2})$ on
$C_3$, the quartic equation of $C_3$ must be:
$$
C_3:\quad Q_1Q_3-Q_2^2=0.
$$
The curve $C_5$
defined by:
$$
s^2=Q_1,\quad st=Q_2,\quad t^2=Q_3 \qquad(\subset\PP^4)
$$
is a canonically embedded genus 5
curve in $\PP^4$ (with coordinates $x,y,z,s,t$ and where $x,y,z$
are a basis of $H^0(C_3,K_3)$).

Projection onto $\PP^2$ defines a 2:1 unramified covering $\pi:C_5\rightarrow
C_3$ and clearly $\pi^*( K_{3}\otimes\alpha)\cong  K_{C_5}$, so $\pi$ is
defined by $\alpha$ and $C_5\cong D_E$.

 From the theory of Prym varieties we have (cf.\ \cite{Mumford}, $\S$ 6):
$$
Nm^{-1}(K_3)=P^+\cup P^-
\qquad{\rm with}\quad Nm:\; Pic^4(C_5)\longrightarrow Pic^4(C_3)
$$
and $P^+,\,P^-$ are both isomorphic to $J(C)$, the Prym variety of the cover
$C_5\rightarrow C_3$ (\cite{Verra}, p.\ 438).
Here we have:
$$
P^+:=\{ L\in Pic^4(C_5):\; Nm(L)=K_3,\;\;
\dim H^0(C_5,L)\equiv 0\;{\rm mod}\;2\},\quad{\rm and}\quad
\tilde{\T}\cap P^+=\Xi,
$$
where $\tilde{\T}$ is the theta divisor in $Pic^4(C_5)$ and where $\Xi$ is the
theta divisor of the Prym variety (actually the intersection has multiplicity
2), so in our case $\Xi=C$.

A point of $C$ thus corresponds to a $g^1_4$ on $C_5$ with norm $K_3$. The
$g^1_4$'s on $C_5$ are cut out by rulings of quadrics in the ideal of $C_5$ of
rank $\leq 4$. The hyperelliptic involution on $C$ corresponds to the
permutation
of the rulings in the rank 4 quadrics, so the Weierstrass points correspond to
the $g^1_4$'s from rank 3 quadrics.

To a section $as+bt\in H^0(K_3\otimes\alpha)$ corresponds
a quadric of rank $\leq 4$ in the ideal of $C_5$ given by:
$$
(as+bt)^2=a^2Q^2+2abQ_2+b^2Q_3,\quad{\rm so}\quad
(as+bt)^2=Q_{(a:b)}.
$$
Such a quadric has rank 3 iff $det(Q_{(a:b)})=0$. Thus these rank 4 quadrics
are parametrized by $|K_3\otimes\alpha|$ and there are 6 rank 3 quadrics that
correspond to the pairs of bitangents.
Each quadric is a cone over a 2:1 cover of the plane $s=t=0$ branched along the
conic $Q_{(a:b)}$. The rulings of a rank four quadric
are the two irreducible components in the inverse image of lines tangent to the
conic $Q_{(a:b)}$; they are interchanged by the covering involution
($s,\,t\mapsto -s,\,-t$).

Any $\PP^2$ in such a rank 4 quadric thus projects to a line tangent to the
conic   $Q_{(a:b)}$ and the divisor cut out by the $\PP^2$ on $C_5$ maps onto
the divisor cut out on $C_3$ by that tangent line.
Hence the norm of the $g^1_4$'s obtained from these quadrics is $K_3$.
This shows that $C$ is indeed the double cover of $|K_3\otimes\alpha|$ branched
over the
six points corresponding to pairs of bitangents.
\qed

\subsection{}
The proposition allows us to make `consistent' choices for the coordinate
functions of the Hitchin map as $E$ varies.
Let $s_i\in H^0(C_3,K_3\otimes\al)$ be the six sections which correspond to
bitangents. Then  $H_i$ restricts to $s_i^2$ on $C_3$,
so if we put:
$$
\bar{H}:\PP T^*_E{\cal M} \longrightarrow \PP^2,\qquad
p\mapsto (H_1(p):H_2(p):H_3(p))
$$
then we have a choice of coordinate functions for $\bar{H}$ which makes sense
for any (general) $E$.

The only remaining problem is that we can can still multiply each $H_i$ by a
function on ${\cal M}\cong\PP^3$ which has poles and zeros in the locus where
the map $D_E\rightarrow \bar{D}_E$ is not an
\' etale 2:1 map of smooth curves.

\section{Computing the Hitchin map}\label{compute}

\subsection{}
In the previous section we saw that the
polynomials $H_i$ on $\PP^3\times (\PP^3)^*$ defining the Hitchin map
have the property:
for any  general $q\in\PP^3$,
$$
(H_i=0)\cap \PP T^*_q\PP^3=l_i\cup l_i'\qquad
(\PP T^*_q\PP^3\subset\PP T^*\PP^3=I=\{(x,h)\in\PP^3\times(\PP^3)^*:\;x\in
h\;\}),
$$
where $l_i$ and $l_i'$ form the pair of bitangents to the smooth curve
$$
C_3:=S\cap \PP T^*_q\PP^3
$$
(see \ref{kum}) corresponding to $s_i\in H^0(C_3,K_3\otimes\alpha)$
(i.e.\ $(s_i^2)=C_3\cap (l_i\cup l_i'$)). Here $S\subset (\PP^3)^*$
is the Kummer surface of $Pic^{g-1}(C)$ and $\alpha$ is the bundle of order two
defined by the \' etale double cover of $C_3$ obtained by pull-back from the
map
$Pic^{g-1}(C)\rightarrow S$.

We now consider the problem of finding such polynomials.

\subsection{}
This problem was actually solved a century ago using the relation between
Kummer surfaces and Quadratic line complexes. The classical solution
is as follows.

A line in $(\PP^3)^*$ with Klein coordinates (see \ref{linec})
$(x_1:\ldots:x_6)\in\PP^5$ is a bitangent to the Kummer surface $S$
occurring in one of the six pairs
iff there is an $i\in\{1,\ldots ,6\}$
such that the following two equations are satisfied (see Proposition
\ref{bip}):
$$
x_i=0,\qquad
\sum_{j\neq i} \frac{x_j^2}{\la_i-\la_j}=0
\qquad\qquad(j\in\{1,\ldots ,6\}),
$$
where the $\la_i$ correspond to the Weierstrass points of the curve $C$:
$$
C:\quad y^2=(x-\la_1)\ldots (x-\la_6).
$$
We will show in \ref{biteq} how to derive the $H_i$ by `restricting'
these two equations to the incidence bundle $I$.

\subsection{}\label{linec} We start with some definitions from line geometry.
The Pl\"ucker coordinates of the line $l=\langle
(Z_0:\ldots: Z_3),\;(W_0:\ldots :W_3)\rangle\subset(\PP^3)^*$ are:
$$
p_{ij}:=Z_iW_j-W_iZ_j\quad{\rm and}\quad
G:\;\;p_{01}p_{23}-p_{02}p_{13}+p_{03}p_{12}=0
$$
is the equation of the Grassmannian of lines, embedded in $\PP^5$.
The Klein coordinates of a line are:
$$
\begin{array}{lll}
X_1=p_{01}+p_{23},\quad& X_3=i(p_{02}+p_{13}),\quad&X_5=p_{03}+p_{12}\\
X_2=i(p_{01}-p_{23}),&X_4=p_{02}-p_{13},&X_6=i(p_{03}-p_{12}).
\end{array}
$$
Note that each $X_i$ corresponds to a non-degenerate alternating bilinear form
in the $Z_i,\,W_i$. These six bilinear forms give sections $\Phi_i$ of the
bundle projection $\PP T^*\PP^3\longrightarrow \PP^3$:
$$
\Phi_i:\PP^3\longrightarrow \PP T^*\PP^3=I\subset\PP^3\times(\PP^3)^*,
\quad q\mapsto (q,\e_i(q)):=(q,X_i(q,-)).
$$
That $\Phi_i(q)\in I$ follows from the fact that $X_i$ is alternating:
$X_i(q,q)=0$. Explicitly, if $q=(x:y:z:t)\in\PP^3$, then the
$\e_i=\e_i(q)\in\PP^{3*}$ have the dual coordinates:
$$
\begin{array}{lll}
\e_1=(y:-x:t:-z),\quad& \e_3=(z:t:-x:-y),\quad&\e_5=(t:z:-y:-x)\\
\e_2=(y:-x:-t:z),&\e_4=(z:-t:-x:y),&\e_6=(t:-z:y:-x).
\end{array}
$$

\subsection{}\label{biteq}
We show how to take care of the first equation.
Let $q=(x:y:z:t)\in\PP^3$.
As the incidence bundle is the cotangent bundle we have:
$$
T^*_q\PP^3=\{(u:v:w:s)\in (\PP^3)^*:\; xu+yv+zw+ts=0\;\}.
$$
Note that we can rewrite the equation to obtain:
$$
T^*_q\PP^3=\{p\in(\PP^3)^*:\;X_i(\e_i(q),p)=0\;\}.
$$
This implies that the lines in $T^*_q\PP^3$ with $X_i=0$ (which form a linear
line complex) are exactly the lines passing through the point $\e_i(q)$
(cf.\ \cite{GH}, p. 759-760).

In particular, if $T^*_q\PP^3\cap S$ is a smooth quartic curve,
and $l_i,\;l'_i$ are a pair of bitangents as before,
then both lines have $X_i=0$ and thus they must intersect in $\e_i(q)$.

Let now $p\in T^*_q\PP^3,\; p\neq \e_i(q)$.
The condition that $p\in l_i\cup l'_i$ is  equivalent to demanding that the
line $\langle\e_i(q),p\rangle$
is one of these two bitangents.
The $i$-th Klein coordinate of this line is zero because it passes through
$\e_i(q)$.

Thus for these lines the first equation is verified. We conclude:
$$
p\in l_i\cup l_i'\subset \PP T^*_q\PP^3\;
\Longleftrightarrow\;
H_i(p,q):=\sum_{j\neq i} \frac{x_j^2}{\la_i-\la_j}=0,\quad
{\rm with}\quad x_j:=X_j(\langle\e_i(q),p\rangle),
$$
the Klein coordinates of the line $\langle\e_i(q),p\rangle\subset (\PP^3)^*$.

The coordinates of
$\e_i(q)$ are linear in those of $q$ and so the Pl\"ucker coordinates of
$\langle\e_i(q),p\rangle$
are homogeneous of bidegree (1,1) in those of $q$ and $p$. Thus $H_i$ is given
by a homogeneous polynomial of bidegree (2,2) on $\PP^3\times (\PP^3)^*$.

\subsection{}
On the open subset $T^*\CC^3=\CC^3\times (\CC^3)^*$ of $T^*\PP^3$
one can obtain a CIHS from the polynomials $H_i$ as follows.

Let $(x,y,z)$ be coordinates on $\CC^3$ and let $(u,v,w)$ be the dual
coordinates on $\CC^{3*}$. Then, the inclusion of cotangent bundles
followed by the (rational) projectivization map is given by
$$
T^*\CC^3 \longrightarrow T^*\PP^3\longrightarrow \PP T^*\PP^3
$$
$$
(q,p):=((x,y,z),(u,v,w))\mapsto
(\tilde{q},\tilde{p}):=((x:y:z:1),(u:v:w:-(xu+yv+zw))).
$$
(Note that the last coordinate is obtained from the incidence condition.)
Now we define:
$$
H^a_i(p,q)=\sum_{j\neq i}
\frac{X_j(\langle
\e_i(\tilde{q})\,,\,(u:v:w:-(xu+yv+zw))\rangle)^2}{\la_i-\la_j}.
$$
The $H_i$ are homogeneous of degree (2,2), and the last coordinate has degree
one in $x,\,y,\,z$ so the $H^a_i$ will have degree $\leq 4$ in the $x,\,y,\,z$
(and need not be homogeneous in these variables), but they
are still homogeneous of degree $2$ in the $u,\,v,\,w$.

With the help of a computer,
one can explicitly write down the polynomials $H^a_i$ (the expressions are
rather long though).
To verify that these polynomials actually Poisson commute (with respect to the
standard two form
$dx\wedge du+dy\wedge dv+dz\wedge dw$) we again used the computer (after
normalizing
 three of the $\la_i$'s by a linear fractional transformation). This then
allows us to
conclude that the map $H^a:T^*\CC^3\rightarrow \CC^3$ (whose coordinate
functions are any three of the six $H_i^a$'s) is a CIHS.

It seems reasonable to expect that the CIHS defined by these $H^a_i$
is actually Hitchin's system, but we could not establish that.

\section{Quadratic Line Complexes}\label{Qcomp}

\subsection{}
In this section we recall how the equations for the bitangents are determined.
We summarize the results we need from \cite{GH}, Chapter 6 and follow
\cite{Hu}.

Let $G\subset \PP^5$ be the Grassmannian of lines in $\PP^3$, so $G$ is
viewed as a quadric in $\PP^5$. For $x\in G$
we denote by $l_x$ the corresponding line in $\PP^3$. For $p\in\PP^3$  and
$h\subset \PP^3$ a plane we
define
$$
\s(p):=\{x\in G:\; p\in l_x\},\qquad \s(h):=\{x\in G:\; l_x\subset h\}.
$$
Both $\s(p)$ and $\s(h)$ are isomorphic to $\PP^2$, in fact any (linear)
$\PP^2$ in $G$ is either a $\s(p)$ or a $\s(h)$. Let $L$ be a line
in $G$, then $L$ is the intersection of a (unique) $\s(p)$ with
a (unique) $\s(h)$:
$$
L=\s(p)\cap\s(h)=\{x\in G:\; p\in l_x\subset h\}.
$$
We will sometimes write $h=h_L$, $p=p_L$ and $L=L_{p,h}$.
Thus the points on the line $L$ (in $G$) correspond to the lines (in $\PP^3$)
in a pencil in
$h$ with `focus' $p$.

\subsection{} A quadratic line complex $X$ is the intersection of $G$ with
another quadric $F$; we assume $X$ to be smooth.
$$
X:=G\cap F.
$$
For any $p\in\PP^3$,
the intersection of $\s(p)=\PP^2\subset G$ with the quadric $F$ is a conic
in
$\s(p)$. Let
$$
S:=\{p\in\PP^3:\; \s(p)\cap F\;{\rm is}\;{\rm singular}\},
$$
then $S$ is a Kummer surface.

\subsection{}
If $\s(p)\cap F$ is singular, it is the union of two
lines $L,\,L'$ or it is a double line. The double lines correspond to the 16
singular
points of $S$. The points $x\in L$ correspond to the lines in a plane
$h_L$ passing through $p$, similarly the points in $L'$ correspond to lines
in a plane $h_{L'}$ passing through $p$. These pencils are called
confocal pencils (having the same focus $p$). Note that the line $l=h_L\cap
h_{L'}$ lies in both these pencils; it is called a singular line of the
complex $X$. This line $l$ corresponds to the intersection of $L$ and
$L'$ in $G$: $[l]=L\cap L'$.

 The singular lines of $X$ form a smooth surface $\Sigma$ in $G$.
$$
\Sigma:=\{x\in X:\;l_x\;\mbox{is a singular line in}\;X\}.
$$
The set $\Sigma$ is determined in \cite{GH}, p.\ 767-769:
$$
x\in\Sigma\; \Longleftrightarrow\; T_xF=T_{x'}G\quad \mbox{for some}\;x'\in G.
$$

\subsection{}\label{kc}
In Klein coordinates $X_i$ the relation between the points $x$ and $x'$ above
assumes a very simple form. Any quadratic line complex $X$ can be given by
(\cite{GH},p.\ 789):
$$
G:\;\; X_1^2+X_2^2+\dots +X_6^2=0,\quad
F:\;\;\la_1X_1^2+\la_2X_2^2+\ldots +\la_6X_6^2=0,
\quad X=G\cap F.
$$
Then $S$ is the Kummer variety associated with the genus two curve
$$
C:\quad y^2=(x-\la_1)\ldots (x-\la_6).
$$
The equations defining the surface $\Sigma$ are then
(\cite{GH}, p.\ 769)
$$
\Sigma=G\cap F\cap F_2,\qquad
F_2:\quad \lambda_1^2X_1^2+\dots +\lambda_6^2X_6^2=0.
$$
Let now
$$
x=(x_1:\ldots:x_6)\in \Sigma\;\Longrightarrow\;\;
T_xF:\quad\lambda_1x_1X_1+\ldots+\lambda_6x_6X_6=0.
$$
Defining
$$
x':=(\lambda_1x_1:\ldots:\lambda_6x_6),\;\Longrightarrow\;\;
\quad
x'\in G,\quad
T_{x'}G:\quad\lambda_1x_1X_1+\ldots+\lambda_6x_6X_6=0,
$$
so $T_xF=T_{x'}G$ and $x'$ satisfies the required condition.

\subsection{Lemma.}\label{lembit}
Let $x\in \Sigma\subset G$ and let $x'\in G$ as above. Define a line:
$$
L:=\langle x,\,x'\rangle\;\subset \PP^5.
$$
Then we have $L\subset G$ and
$$
L=L_{p,h}\qquad {\rm with}\quad p\in S,\quad h=T_pS,
$$
so that the points $y\in L$ correspond to the lines $l_y\subset \PP^3$ with
$p\in
l_y\subset T_pS$. For $i\in\{1,\ldots ,6\}$ let
$$
\{[l_i]\}:=\;(X_i=0) \cap L\qquad(\in G\subset \PP^5).
$$
Then $l_i$ is a bitangent line to $S$. Moreover,
if $p$ is a general point of $S$ then any bitangent to $S$ passing
through $p$ is one of the six $l_i$'s.

\ts
Since $x\in T_{x'}G$ we have $L\subset G$ (this is also easily verified using
the three equations defining $\Sigma$). Then $L=L_{p,h}$ with $p=l_x\cap
l_{x'}$.
We claim that $p\in S$. Since $\s(p)$ is a linear subspace in $G$ passing
through $x'$ we have $\sigma(p)\subset T_{x'}G$, and thus also
$\sigma(p)\subset T_xF(=T_{x'}G)$. Thus $\s(p)$ is tangent to $F$ at $x$, so
$\s(p)\cap F$ is singular in $x$. Therefore $p\in S$ (and $l_x$ is the singular
line of $X$ passing through $p$).

Any point on $L$, distinct from $x$, can be written as:
$$
x_\la:=\la x+x'=(\ldots:(\la +\la_i)x_i:\ldots)\qquad (\la\in {\CC}).
$$
It is easy to check by substitution that $x\in \Sigma\Rightarrow
x_\la\in\Sigma_\la$ with:
$$
\Sigma_\la:= G\cap F^{(\la)}\cap F^{(\la)}_2,\quad(\la\neq -\la_i)
$$
and where we define:
$$
F^{(\la)}:\;(\la+\la_1)^{-1}X_1^2+\ldots+(\la+\la_6)^{-1}X_6^2=0,\quad
F^{(\la)}_2:\;(\la+\la_1)^{-2}X_1^2+\ldots+(\la+\la_6)^{-2}X_6^2=0.
$$
Thus $x_\la$ corresponds to a singular line for the quadratic complex
$X_\la:=G\cap F^{(\la)}$.

As above, there exists thus a point $x'_\la\in G$ with:
$$
T_{x_\la}F^{(\la)}=T_{x'_\la}G,\qquad
x'_\la:=(\ldots:x_i:\ldots)=x\in G.
$$
Therefore $x_\la$ and $x'_\la$ lie on the line $L\subset G$ and thus the point
$p$ is a point of $S_\la$, the Kummer surface associated to the quadratic line
complex $X_\la$. This holds for all singular lines $x$ of $X$ (and thus for all
points $p\in S$), therefore we conclude:
$$
S=S_\la.
$$
In particular, there is a one-dimensional family of quadratic line complexes
$X_\la$ which give rise to the same Kummer surface, the so called Klein variety
(see \cite{NR} for a modern treatment).

Now we can determine $h$. Each $x_\lambda\in L$ is a singular line for a
quadratic line complex defining $S$. Then the line in $\PP^3$ corresponding to
it is tangent to $S$ at the (unique; cf. the verification on p.\ 767 of
\cite{GH})
point $p_\lambda\in S$
with $x_\la=Sing(\sigma(p_\la)\cap F)$ (cf.\ \cite{GH},
p.\ 764-765, p. 791).
In our case, $p_\la=p$ for all $\la$, so we conclude that $L$ is the pencil of
lines in $\PP^3$ that are tangent to $S$ at $p$, which implies $h=T_pS$.

The lines from $L$ that are bitangent to $S$ are thus the bitangents of the
curve $T_pS\cap S$.
This is, in general, a plane quartic curve with a node at $p$, so there are six
bitangents to $S$ in the pencil $L$. These must then correspond to the values
$\la=-\la_i$, since for other values the lines $x_\la$ are singular lines of a
smooth quadratic line complex and cannot be tangent to $S$ at other points.
Thus the bitangent lines in $T_pS$ passing through $p$
correspond to the points on $L$ with exactly one Klein coordinate
equal to zero.
\qed

\subsection{Remarks.} Viewing $S$ as $Pic^{g-1}(C)/\iota$, the divisors
$T_pS\cap S$ (for $p$ smooth) correspond to the divisors
$$
D_\beta:=L_\beta^*\T+L_{-\beta}^*\T\in |2\T|,\qquad(
\beta \in Pic^0(C))
$$
with $2\beta\neq\cO$.
These are the union of two copies of $C\,(=\T)$ meeting in two points. These
two points, and the two copies, are interchanged by $\iota$, the quotient is a
nodal curve isomorphic to $C$ with the two points identified.
The normalization of $T_pS\cap S$ is thus isomorphic with $C$. The six
bitangents to $C$ correspond to the lines spanned by $p$ and (the image in
$T_pS\cap S$) of a Weierstrass point of $C$.
This, once again, establishes the connection between bitangents to $S$
(and its plane sections) and Weierstrass points on $C$.

The 16 non-reduced divisors $L_\beta^*(2\T)\in |2\T|$ (so $2\beta=\cO$) map to
double conics.
A point on a double conic is `exceptional' for the Lemma since
any line tangent to $S$ at a point of the conic lies in the plane of the conic
and is thus a bitangent to $S$.
These (double) conics are called the tropes
of the Kummer surface.

The `self'-duality of the Kummer surface $S=\delta(Pic^{g-1}(C))$
fits in nicely with the map
$\delta:Pic^{g-1}(C)\rightarrow |2\T|_0$ from \ref{delta}, the duality
between $|2\T_0|$ and $|2\T|$, and the map $\delta':Pic^0(C)\rightarrow
|2\T|,\;
\beta\mapsto D_\beta$. In fact, it identifies the tangent planes to points of
$S$ (which cut out $D_\beta$) with the points $\pm \beta \in Pic^0(C)/\pm 1
\cong \delta'(Pic^0(C))$ which is the Kummer surface of $Pic^0(C)$.
This surface
is isomorphic to $S=Pic^{g-1}(C)/\iota$, but the `self'-duality
is however not an
isomorphism; it is a birational isomorphism which blows up
double points and blows down tropes.

The special case that the line $L$ in the Lemma
actually lies in $X$ is studied in \cite{GH}, p.\ 791-796 (note that they fix
the
quadratic line complex whereas in the proof of the Lemma we consider a family
of complexes).

\subsection{Proposition.}\label{bip}
For $x\in G$ the line $l_x$ is a bitangent to a general point of $S$
iff for some $i\in\{1,\ldots ,6\}$ one has:
$$
x_i=0,\qquad
\sum_{j\neq i} \frac{x_j^2}{\la_i-\la_j} =0.
$$

\vspace{\baselineskip}
\noindent{\bf Proof.}$\;\;$
In Lemma \ref{lembit} we saw that any such
bitangent $l_z$ of $S$ has one Klein coordinate equal to zero,
we will assume it is the first one. Then
$$
z=(0:(-\la_1+\la_2)x_2:\ldots :(-\la_1+\la_6)x_6)\quad
{\rm with}\quad
x=(x_1:x_2:\ldots:x_6)\in\Sigma,
$$
in particular $x\in X=G\cap F$.
Substituting the coordinates of $z$ in the second equation we get:
$$
(\la_1-\la_2)x^2_2+\ldots +(\la_1-\la_6)x_6^2=\la_1(x_1^2+\ldots +x_6^2)-
(\la_1x_1^2+\ldots+\la_6x_6^2)=0.
$$

Conversely, let $z=(0:z_2:\ldots :z_6)\in G$ satisfy also the second equation
above, so:
$$
z_2^2+\ldots +z_6^2=0,\quad
(\la_1-\la_2)^{-1}z^2_2+\ldots +(\la_1-\la_6)^{-1}z_6^2=0.
$$
Then we define:
$$
x_i:=(\la_1-\la_2)^{-1}z_i,\quad(2\leq i\leq 6),\qquad
x_1:=\sqrt{-(x_2^2+\ldots +x_6^2)},
$$
here the choice of the square root does not matter. Define
$$
x:=(x_1:x_2:\ldots :x_6),\quad{\rm so}\quad x_1^2+\ldots +x_6^2=0
$$
and we have $x\in G$ (the quadric defined by $X_1^2+\ldots +X_2^2$.)

We claim that $l_x$ is a singular line of $S$. For this we verify that $x$
satisfies the other two quadratic equations defining $\Sigma$. First of all:
$$
\begin{array}{rcl}
\la_1x_1^2+\la_2x_2^2+\ldots+\la_6x_6^2&=&
-\la_1(x_2^2+\ldots +x_6^2)+\la_2x_2^2+\ldots+\la_6x_6^2\\
&=&(\la_1-\la_2)x_2^2+\ldots+(\la_1-\la_6)x_6^2\\
&=&(\la_1-\la_2)^{-1}z^2_2+\ldots
+(\la_1-\la_6)^{-1}z_6^2\\
&=&0,
\end{array}
$$
so $x\in F$. We use these two relations on the $x_i$'s to obtain the third:
$$
\begin{array}{rcl}
\la_1^2x_1^2+\ldots+\la_6^2x_6^2&=&
(\la_1^2x_1^2+\ldots+\la_6^2x_6^2)-2\la_1(\la_1x_1^2+\ldots+\la_6x_6^2)
+\la_1^2(x_1^2+\ldots+x_6^2)\\
&=&
(\la_1-\la_2)^2x_2^2+\ldots +(\la_1-\la_6)^2x_6^2\\
&=& z_2^2+\ldots +z_6^2\\
&=&0,
\end{array}
$$
since $z=(0:z_2:\ldots:z_6)\in G$. Thus $x\in F_2$ and we conclude
$x\in\Sigma$,
so we verified that $l_x$ is a singular line.
Note that (with notation from \ref{kc}):
$$
x'=(\la_1x_1:\ldots:\la_6x_6),\quad {\rm so}\;\; z=-\la_1x+x'
$$
and thus $l_z$ is indeed a bitangent to $S$ (see the proof of the Lemma).
\qed

\

\

\

\newbox\ik
\newbox\zij
\setbox\ik\hbox{email: geemen@math.ruu.nl }
\setbox\zij\hbox{Department of Mathematics}

\parbox{\wd\ik}
{
Bert van Geemen \\
Department of Mathematics \\
P.O.\ Box 80.010 \\
3508 TA Utrecht \\
The Netherlands \\
\mbox{}  \\
e-mail: geemen@math.ruu.nl
}
\hfill
\parbox{\wd\zij}
{Emma Previato \\
Department of Mathematics \\
Boston University \\
Boston, MA 02215 \\
USA \\
\mbox{}  \\
e-mail: ep@math.bu.edu
}


\begin{thebibliography}{ACGH}

\bibitem[B1]{B1} A. Beauville, Fibr\'es de rang 2 sur une courbe, fibr\'e
d\'eterminant et fonctions th\^eta, {\it Bull. Soc. Math. France} {\bf
116} (1988), 431-448.

\bibitem[B2]{B2} A. Beauville, ---, II,{\it ibid.} {\bf 119} (1991), 259-291.


\bibitem[Be]{B} A. Bertram, Moduli of rank-2 vector bundles, theta
divisors, and the geometry of curves in projective space, {\it J.
Differential Geom.} {\bf 35} (1992), 429-469.

\bibitem[BNR]{BNR} A. Beauville, M.S. Narasimhan and S. Ramanan, Spectral
curves and the generalized theta divisor, {\it J. reine angew. Math.}
{\bf 398} (1989), 169-179.

\bibitem[F]{F}  G. Faltings, Stable G-bundles and projective connections,
{\it J. Alg. Geometry} {\bf 2} (1993), 507-568.

\bibitem[GH]{GH} P. Griffiths and J. Harris, {\it Principles of Algebraic
Geometry}, Wiley, 1978.

\bibitem[H1]{H1} N.J. Hitchin, The self-duality equations on a Riemann
surface, {\it Proc. London Math. Soc.} {\bf 55} (1987), 59-126.

\bibitem[H2]{H2} N.J. Hitchin, Stable bundles and integrable systems,
{\it Duke Math. J.} {\bf 54} (1987), 91-114.

\bibitem[Hu]{Hu} R.W.H.T. Hudson, {\it Kummer's quartic surface}, Cambridge
Univ. Press, 1990.

\bibitem[L]{L} Y.  Laszlo, Un th\' eor\` eme de Riemann pour les diviseurs
th\^ eta sur les espaces des modules de fibr\' es stables sur une courbe,
{\it Duke Math. J.} {\bf 64} (1991), 333-347.

\bibitem[M]{Mumford} D.\ Mumford,  Prym varieties I, in {\it Contributions
to Analysis}, Academic Press, 1974, pp. 325-350.

\bibitem[NR1]{NR} M.S. Narasimhan and S. Ramanan, Moduli of vector bundles
on a compact Riemann surface, {\it Ann. of Math.} {\bf 89} (1969),
19-51.


\bibitem[NR2]{NR2} M.S. Narasimhan and S. Ramanan,
2$\T$ linear systems on Abelian varieties, in {\it Vector Bundles
on algebraic varieties}, Oxford Univ. Press, 1987, pp. 415-427.

\bibitem[V]{Verra} A. Verra, The fibre of the Prym map in genus three,
{\it Math. Ann.} {\bf 276} (1987), 433-448.


\end{thebibliography}
\end{document}